\documentclass{amsart}
\usepackage{amssymb}
\newtheorem{theorem}{Theorem}
\newtheorem{lemma}{Lemma}
\newtheorem{definition}{Definition}
\def\mR{\mathcal R}
\begin{document}
$\, $
\bigskip
\bigskip
\bigskip
\bigskip

\begin{center}{\Large Polynomial birth--death processes and the second conjecture of Valent}
\end{center}
\bigskip \bigskip

\hskip2,5cm\vbox{\hsize10,5cm\baselineskip4mm
\noindent{\small\textbf{Abstract.} The conjecture of Valent about the type of Jacobi matrices with polynomially growing weights is proved.}

\noindent\textbf{Keywords:} moment problems, birth--death processes,  Jacobi matrices.}

\footnotetext{AMS subject classifications: 34L15, 47B36.}

\bigskip

\centerline{\sf Ivan Bochkov}

\medskip

\begin{center}
Faculty of Mathematics and Mechanics, \\
St Petersburg State University, \\
198504, St Petersburg, Russia,\\
e-mail: bv1997@ya.ru\\
\end{center}

\bigskip

\section{Introduction} 
In 1998 G. Valent in \cite{Valent} conjectured the order and type of certain
indeterminate Stieltjes moment problems associated with birth and death
processes having polynomial birth and death rates of degree $p \ge 3 $. His conjecture says that the order of the birth--death processes with rates $ \lambda_n $, $ \mu_n $ being the polynomials of degree $ p $, 
\[ \lambda_n = n^p + C n^{ p-1 } + \dots , \] \[ \mu_n = n^p + D n^{ p-1 } + \dots ,\] subject to the condition $ 1 < C - D < p-1 $ is $ 1/p $, and its type, $ t_p $, with respect to that order is 
\begin{equation}\label{valentconj} p \int_0^1 ( 1 - x^p )^{ - 2/p } d x , 
\end{equation} 
respectively. The conjecture was formulated on the basis of several explicitly solvable examples for $ p=3 $ and $ p=4 $ found by Valent and his collaborators, see \cite{BergValent}, \cite{GLV}.

In \cite{BergSwarc} the conjecture has been reduced to the following question in terms of Jacobi matrices. Let $ a_j = j^p $, $ p> 1 $, $ j \ge 1 $, and let
\[ J_p = \begin{pmatrix} 0 & a_1 & 0 & \cdots & \cdots \cr
a_1 & 0 & a_2 & 0 & \cdots  \cr
0 & a_2 & 0 & a_3 & \cdots \cr 
\cdots & 0 & a_3  & 0 & \cdots \cr
\cdots & \cdots & \cdots & \cdots & \cdots \end{pmatrix} . \]
Then the corresponding Hamburger moment problem is well-known to be indeterminate \cite[Chapter 1, Exercise 5]{Akhiezer}. As shown in \cite{BergSwarc}  the Valent conjecture holds true if the elements of the corresponding Nevanlinna matrix have order $ p^{ -1 } $ (first conjecture) and the type with respect to that order is $ 2^{ -1 } B ( 1/(2p) , 1 - 1/p ) $ (second conjecture). Here and throughout the paper $ B $ is the Euler beta-function.

The conjecture about the order was proved by Romanov in \cite{Rom} as an application of a general method of estimating the order of canonical systems developed in that paper. In \cite{BergSwarc} Berg and Szwarc gave another proof of the order conjecture and established that the type satisfies  
\[ \frac\pi{ \sin ( \pi / p ) } \le t_p \le \frac\pi{ \sin ( \pi / p ) \cos ( \pi / p )  } . \]
This estimate  is compatible with the second Valent conjecture in the sense that the quantity \eqref{valentconj} satisfies this inequality.

In the present paper the second Valent conjecture is proved completely. Our proof uses the following assertion due to Berg and Szwarc (a minor misprint in the formulation is corrected)

\begin{theorem}[\cite{BergSwarc}, Theorem 1.11]\label{BS}
\begin{equation} \label{BBS} t_p = \frac p{e} \limsup_{ n \to \infty } \left( n  \left(  \sum_{1 \le x_1 \le x_2 < x_3 \le x_4< \ldots < x_{ 2n-1} \le x_{2n}}(x_1x_2\cdots x_{2n})^{-p} \right)^{ \frac1{2pn} } \right).
\end{equation}
\end{theorem}

In the following theorem the sign $ \prec $ in the summation index means $ < $ or $ \le $ depending on the parity of the number $ n $ involved.

\begin{theorem}
Let $p>1$ be a real number, and let 
\[ s(n) = \sum_{1 \le x_1 \le x_2 < x_3 \le x_4< \ldots \prec x_{n}}(x_1x_2\cdots x_{n})^{-p} , \]
  \[ k_n=n (s(n))^{\frac{1}{np}}. \]
Then \[ k_n \mathop{\longrightarrow}_{ n \to \infty } \frac{e}{p}B \left(\frac{1}{2p}, 1- \frac{1}{p} \right) . \]
\end{theorem}

As follows from Theorem \ref{BS} above this assertion implies the Valent conjecture. Before procceding to the proof of Theorem 2, let us explain where (\ref{BBS}) comes from, referring to \cite{BergSwarc} for details. One can associate with the Jacobi matrix $ J_p $ a sequence of $ 2  \times 2 $ elementary monodromy matrices, $ M_j $,  of the form $ M_j = I + z R_j $ where $ z $ is a complex parameter and $ R_j $ is an upper or lower-triagonal matrix, depending on parity of $ j $, explicitly calculated in terms of $ a_j $. The infinite product $\dots M_n \dots M_2 M_1 $ is essentially the Nevanlinna matrix. On developing this product of sums and taking into account the triangle structure we end up with  
an explicit expression for Taylor coefficients (in $ z $) of the matrix elements. Theorem \ref{BS}  is just an expression for the type of those functions via the Taylor coefficients.

Let us mention a wide context of the Valent conjecture. It belongs to the theory of indeterminate moment problems \cite{Akhiezer}. According to the M. Riesz theorem, the Nevanlinna matrices corresponding to indeterminate problems have minimal exponential type (with respect to the order $1$). This leads to the question about their actual order and type with respect to that order. The examples where the order and especially the type are known are few and isolated. Apart from those already mentioned most of them come from explicitly solvable orthogonal polynomials within the $q$--Askey scheme \cite{Koekoek} and have order zero. The main difficulty is high instability of the indeterminate problems. In particular, the estimates obtained by the variational approach (minimaximal principle) to the spectrum of the corrresponding Jacobi matrices are not precise enough to calculate the type. 

\section{Proof of Theorem 2}
 
\subsection{Preliminaries} From now on, we fix the number $p$ and do not indicate the dependence of irrelevant constants from it. 

\begin{definition} Two sequences,  $x_n$ and $y_n$ of positive reals are said to be equal, denoting $ x_n \approx y_n $, if $ \ln ( x_n / y_n ) = o ( n ) $ as $ n \to \infty $. 
\end{definition}

Obviously, if $ x_n \approx y_n $ then the sequences $ n x_n^{ \frac 1{np}} $ and $ n y_n^{ \frac 1{np} }$ converge simultaneously and, if they do, the limits coincide, hence the term.
The structure of the proof of Theorem 1 is as follows. First we describe the steps of the proof and then, in a separate section, provide details for each step. The steps and their details are enumerated accordingly. Before doing so, let us establish for a future reference a ``trivial'' estimate for $ s ( n ) $. 

\begin{lemma}
There exists $ C_1, C_2 > 0 $ such that 
\begin{equation}\label{trivial}  C_1^n n^{ -pn } \le  s( n ) \le C_2^n n^{ -pn }  . \end{equation}
\end{lemma}

Let us first comment on the upper estimate. The summation indices, $ x_j $, in the definition of $ s( n ) $ satisfy  $ x_j \ge \lfloor j/2 \rfloor $, hence $ s ( n ) $ is estimated above by the sum over $ x_j \ge \lfloor j/ 2 \rfloor $, which gives on account of the elementary inequality 
\begin{equation} \label{trivi} \sum_{ m=j+1 }^\infty m^{ -p } \le \frac{ j^{ 1-p }}{ p-1 } \end{equation} 
that  for some $ C $ we have
\[ s ( 2n ) \le C^n (n ! )^{ 1-2p }  . \] 
This is worse than the actual upper estimate in (\ref{trivial}) by the extra $ 1 $ in the exponent in the right hand side. The bulk of the proof of the lemma consists in getting rid of this extra factor. 

\begin{proof}
Observe first that taking $ ( x_1 , x_2 , \dots , x_n ) = ( 1 , 1 , 2 , 2 , \dots ) $ in  the definition of $ s( n ) $ implies $ s( 2n ) \ge \left( n ! \right)^{ -2p} $, and $ s( 2n-1 ) \ge \left( (n-1)! n! \right)^{ -2p} $, and the lower estimate in (\ref{trivial}) follows. 

To establish the upper estimate it is sufficient to prove it for $ n $ even, for 
\begin{equation} \label{evenodd}  s( 2n +1 ) \le \frac{( 2n )^{ 1-p }}{ p-1 } s ( 2n ) \end{equation} 
by the argument mentioned after the formulation of the lemma. From now on we assume that $ n $ is even, hence, in particular, $ x_n $ takes values $ x_{ n-1 } , x_{n-1} + 1 , \dots $. We have,  
\begin{eqnarray*} s(n) \le  \sum_{1 \le x_1 \le x_2 < x_3 \le x_4< \ldots \le x_{ n-2 } < x_{n-1}}(x_1x_2\cdots x_{n-1})^{-p} \left[ x_{n-1}^{-p} + \frac{ x_{ n-1 }^{ 1-p }}{ p-1 } \right] \le \\ \frac 2{p-1}\sum_{1 \le x_1 \le x_2 < x_3 \le x_4< \ldots \le x_{ n-2 } <  x_{n-1}}(x_1x_2\cdots x_{n-2})^{-p} x_{ n-1}^{ 1 -2p} 
\end{eqnarray*}
if $ \lfloor \frac { n-1 }2 \rfloor > p-1 $ since $ x_{ n-1 } \ge \lfloor \frac { n-1 }2 \rfloor $, and therefore $ x_{ n-1 } > p-1 $. Proceeding with the estimate we have 
\begin{eqnarray*}   s( n ) \le \frac 2{2 (p-1)^2}\sum_{1 \le x_1 \le x_2 < x_3 \le x_4< \ldots < x_{n-3} \le x_{ n-2} }(x_1x_2\cdots x_{n-2})^{-p} x_{ n-2}^{ 2 -2p} \le \\ \frac 2{2 (p-1)^2}\sum_{1 \le x_1 \le x_2 < x_3 \le x_4< \ldots \le x_{ n-4 } < x_{n-3}}(x_1x_2\cdots x_{n-3})^{-p}  \left[ x_{n-3}^{2-3p} + \frac{ x_{ n-3 }^{ 3-3p }}{3( p-1) } \right] \end{eqnarray*}
If $ x_{ n-3 }> 3 ( p-1) $ then 
\begin{equation}\label{square}
x_{n-3}^{2-3p} + \frac{ x_{ n-3 }^{ 3-3p }}{3( p-1) } \le 2 x_{n-3}^{3-3p} .
\end{equation}
Plugging this one gets that 
\[ s ( n ) \le  \frac {2 \cdot 2}{2 \cdot 3 (p-1)^3}\sum_{1 \le x_1 \le x_2 < x_3 \le x_4< \ldots \le x_{ n-4 } < x_{n-3}}(x_1x_2\cdots x_{n-4})^{-p}  x_{n-3}^{3-4p} 
\]
if $ \frac n2-2 > 3 ( p-1) $ since $ x_{ n-3 }\ge \lfloor \frac {n-3}2 \rfloor = \frac n2-2 $. Repeating this process we get that if $ \frac{n-j}2  \ge j ( p-1 ) $, $ j $ even, then
\begin{eqnarray*} s( n ) \le \frac {2^{ 1 + j/2 }}{ j! ( p-1)^j } \sum_{1 \le x_1 \le x_2 < x_3 \le x_4< \ldots \prec x_{n-j}}(x_1x_2\cdots x_{n-j-1})^{-p}  x_{n-j }^{j - ( j+1  )p } \le \\ \le \frac {2^{ 1 + j/2 }}{ j! ( p-1)^j } \left( \frac{n - j}2 \right)^{ j ( 1 -p ) } s( n-j ) . \end{eqnarray*} 
Thus, there exists a constant $ C> 0 $ such that 
\begin{equation}\label{sn2}  s( n ) \le  \frac {C^n}{ j^j }  \left( n - j \right)^{ j ( 1 -p ) } s( n-j )  , \; j \le \alpha n  , \; \alpha = \frac 1{2p-1} . \end{equation} 
Here we rook into account the Stirling formula in the denominator. It is convenient for us to define $ s ( t ) $ for non-integer $ t $ by 
$ s ( t ) = s ( \lfloor t \rfloor ) $. Then the estimate \eqref{sn2} holds for non-integer $n $ and $ j $ as well in the stated range of $ j $'s on account of \eqref{evenodd}.  With $ j = \alpha n $ it gives \[ s( n ) \le  \frac {C^n}{ n^{ p \alpha n } }  s( n( 1-\alpha ) ) 
\]
Applying this inequality successively $ \sim \frac{ -\ln n }{ \ln ( 1 - \alpha ) } $ times we find that 
\[ s ( n ) \le \frac { C^{ n + n (1-\alpha)+ n ( 1- \alpha )^2 + \dots } }{ n^{p \alpha n } \left( n ( 1 - \alpha ) \right)^{ ( 1 - \alpha) \alpha p n} \left( n ( 1 - \alpha)^2 \right)^{ ( 1 - \alpha)^2 \alpha p n } \dots } . \] Since $ \alpha < 1 $ the numerator is estimated above by $ C^n $ for an appropriate $ C $, and the denominator equals 
\[ n^{ n \alpha p \sum_{ j = 0 }^{ \frac{ -\ln n }{ \ln ( 1 - \alpha ) }} ( 1 - \alpha )^j } \left( 1 - \alpha \right)^{ \alpha p n ( \sum_{ j \ge 0 } j ( 1 - \alpha )^j   + o ( 1 )) } = n^{ -pn + O ( 1 ) } ( 1 - \alpha )^{n ( 1 + o ( 1 ) ) } . \]
Combining these we get the required upper estimate in \eqref{trivial}. 
\end{proof}

\subsection{Plan of the proof}
Step 1 -- Cutting the tails. 
For $ A $ large enough the sequence 
 $$s'_n=  \sum_{1 \le x_1 \le x_2 < x_3 \le \ldots \prec x_{n} \le T(n)}(x_1x_2\cdots x_{n})^{-p} , \; T(n)=n^A , $$ 
satisfies $ s (n)\approx s'_n $. A careful analysis shows that one can take $ A = \frac{5p}{p-1} $ but the exact value of it does not matter for the proof.
 
Step 2 -- ``dyadisation''.  Given an $\alpha>1$ ($\alpha - 1 $ is to be thought of as a small parameter in what follows) define a function $ P $ by $P(n)= \alpha^k $, $ n \in [ \alpha^k , \alpha^{ k+1 } ) $. Let $T'(n)=\alpha T(n)$ and
\begin{equation}\label{sn} s'_n ( \alpha ) =  \sum_{1 \le x_1 \le x_2 < x_3 \le \ldots \prec x_{n}<T'(n)}(P(x_1)P(x_2)\cdots P(x_{n}))^{-p} , \end{equation} .  
$$ S_n (\alpha)=  \sum_{1 \le x_1 \le x_2 < x_3 \le \ldots \prec x_{n}<T'(n)}(x_1x_2\cdots x_{n})^{-p} , $$
\[ k'_n(\alpha)= n(s'_n(\alpha))^{\frac{1}{np}}, \; K_n(\alpha)= n (S_n(\alpha))^{\frac{1}{np}}. \] 
Since $T'(n)>T(n)$, we have $  S_n (\alpha) \approx s(n) $ for every $\alpha $. Then,  the inequality 
\[ S_n(\alpha)  \le s'_n(\alpha) \le \alpha^{np} S_n(\alpha) 
\]
holds because $\alpha^{-1} \le P(x)/x \le 1$. Thus,
\[ K_n(\alpha) \le k_n^\prime  (\alpha) \le \alpha K_n(\alpha) , \] and the theorem will be established if we show that 
\[ k_n^\prime ( \alpha ) \mathop{\longrightarrow}_{ n \to \infty}
\frac{e}{p}B \left(\frac{1}{2p}, 1- \frac{1}{p} \right) ( 1 + O ( \alpha - 1 ) ) . \]

Step 3 - For a given $\alpha>1$ define $l(n)=\lfloor A \log_{\alpha}n \rfloor$. We will drop the argument $ \alpha $ and the prime sign for notation convenience, writing $ s_n = s_n^\prime ( \alpha ) $ from now on. Each $ P ( x_i ) $ in (\ref{sn})  is one of the numbers $1, \alpha,\ldots, \alpha^{l(n)} $. Let $c_i = \# \{ y \colon P ( y ) =\alpha^i \} $, and let $ H(a_0, a_1, \ldots, a_l)$ be the number of tuples $ ( x_1, x_2 , \dots , x_n ) $ such that $1\le x_1\le x_2<\ldots \prec x_n<T'(n)$ and $ \# \{ j \colon P(x_j)=\alpha^{i} \} = a_i $ for all $i$, $ 0 \le i \le l $. 
Then 
$$ s_n=\sum_{(a_0 , a_1, \ldots ,a_l ) \colon  \sum{a_i}=n}H(a_0, a_1, \ldots , a_l )\alpha^{p(-a_1-2a_2-\ldots - la_l)} . $$ 

Step 4 - calculation of $H(a_0, a_1, \ldots , a_l)$. 
The integers from $ 1$ to $T^\prime ( n ) $ are split into $l+1$ groups of respective sizes  $c_0, c_1, \ldots, c_l$, and according to the definition of $ H $, there are $a_i$ numbers chosen from the $i$-th group. With $a_0, a_1, \ldots, a_l$ fixed the choices from different groups are independent, so $H(a_0, a_1, \ldots , a_l) $ is a product of numbers of ways to choose the $ a_i $ numbers, $ y_j $, such that $y_1 \prec y_2 \ldots \prec y_{a_i}$ from the $i$-th group. Consider $i$-th group.  Without loss of generality one can shift the $y_j$'s so that the shifted numbers become integers from the interval $[ 1, c_i ]$ (replacing $ y_j \to \lfloor y_j - \alpha^i \rfloor + 1$). Let us define the following transformation of the  sequence $ \{ y_j \} $, 
\[ ( y_1, y_2 \ldots ,y_{a_i} ) \to ( y_1+m_1, y_2+m_2, \dots, y_{a_i}+m_{a_i} ) , \] 
\[ 
m_j=\left\{ \begin{array}{cc} \lfloor \frac{j-1}{2}\rfloor , & \sum_{k=0}^{k-1}a_{k} \textrm{ odd} \cr \lceil \frac{j-1}{2}\rceil, & \mathrm{otherwise}. \end{array} \right. \] 
This is a bijective transformation between sequences $ \{ y_j \}_1^{ a_i } $ subject to $ y_1 \prec y_2 \ldots \prec y_{a_i}$ and strictly increasing sequences   $x'_1 < x'_2 <\ldots < x'_{a_i}\le m_{a_i} +c_i$. Thus the number of choices for $ i $-th group is equal to the number of strictly increasing sequences   $x'_1 < x'_2 <\ldots < x'_{a_i} \le w_i+c_i$, $w_i=m_{a_i}$, and the latter number is $ c_i+w_i \choose a_i $. Thus,  
\[ H(a_0, a_1, \ldots, a_l)=\prod^l_{i=0} {c_i+w_i \choose a_i } . \]

Step 5 - Replacing $ w_j $ with $ a_j / 2 $.
Let 
\begin{equation} \label{Hprime}  H'(a_0, a_1, \ldots, a_l)=\prod^l_{i=0}\alpha^{-ia_ip} {c_i+ \frac{a_i}2 \choose a_i } \end{equation}
and  
$$s_n^{ (1) } =\sum_{(a_0 , a_1, \ldots ,a_l ) \colon  \sum{a_i}=n} H'(a_0, a_1, \ldots , a_l) . $$ Then $ s^{ (1)}_n \approx s_n $.

Step 6 -  the sum can be replaced by the maximum of the summand. Let  
 $$ s^{ (2)}_n=\max_{(a_0 , a_1, \ldots ,a_l ) \colon  \sum{a_i}=n}H'(a_0, a_1, \ldots , a_l) . $$
Then $ s_n^{ (2)} \approx s^{ (1)}_n$.

Step 7 - the maximum over an integer hyperplane can be replaced by the maximum over the real one. Notice that the right hand side in the definition  (\ref{Hprime}) rewritten as
$$H'(a_0, a_1, \ldots, a_l)=\prod^l_{i=0}\frac{\alpha^{-ia_ip}}{(c_i+\frac{a_i}{2}+1) B(a_i+1,c_i-\frac{a_i}{2}+1)}$$ makes sense for all real $ a_i $, $0\le a_i \le 2c_i$. Let $ \Pi_n = \{ ( x_0 , \dots , x_l ) \in \mathbb R^{ l+1 } \colon 0 \le x_j \le 2 c_j , \, \sum x_j=n \} $, and
$$s^{ (3)}_n=\max_{ \Pi_n }H'(x_0, x_1, \ldots , x_l) . $$ 
Then $ s_n^{ (3)} \approx s^{ (2)}_n$.

Step 8 - the Stirling formula. Let 
 \begin{equation} \label{s4} s^{(4)}_n=\max_{ \Pi_n } \prod_{i=0}^l \exp \left(-ia_ip\ln \alpha +\left(\frac{a_i}{2}+c_i\right) \ln \left(\frac{a_i}{2}+c_i \right)-a_i\ln a_i -\left(c_i-\frac{a_i}{2}\right)\ln\left(c_i-\frac{a_i}{2}\right) \right) . \end{equation}
 Then $ s_n^{ (4)} \approx s^{ (3)}_n$.

Step 9 - calculating the maximum. Let us denote by $ G $  the argument of the exponent in (\ref{s4}). Considering the conditional extremal problem for $ G $ on the hyperplane $ \sum x_j = n $ we first notice that the derivatives $ \partial G / \partial a_i $ are all equal at a critical point. On calculating,    
\begin{eqnarray*}
\frac{\partial G}{ \partial a_i} =-ip\ln\alpha + \frac{\ln (\frac{a_i}{2}+c_i)}{2}-\ln a_i +\frac{\ln (c_i-\frac{a_i}{2})}{2} = -ip \ln \alpha + \frac 12 \ln \left( \frac{c_i^2}{a_i^2} - \frac 14 \right). \end{eqnarray*} 
The right hand side is a monotone decreasing function of $ a_i $ on $ ( 0 , 2 c_i ) $, going to $ + \infty $ as $ a_i \to 0+ $, to $ -\infty $ as $ a_i \to 2c_i - $ hence it takes any value, $ \lambda $, at exactly one point,
\[ a_i ( \lambda ) =  \frac{2c_i}{\sqrt{4\alpha^{2ip} e^{ 2 \lambda}+1}} . 
\]
It follows that there is exactly one critical point, which is determined from the equation, $\sum{a_i ( \lambda) }=n$, and it is easy to see that this point is the required point of maximum. On pluggung this maximum into (\ref{s4}) we get
\begin{equation}\label{s4res} s^{(4)}_n = \exp \left( \lambda n + \sum^l_{ i=0}  c_i \ln \frac{ \sqrt D_i + 1 }{ \sqrt D_i - 1 } \right) , \end{equation}
\[ D_i \colon = 1 + 4 \alpha^{ 2ip } e^{ 2 \lambda } . \]

Step 10 - The equation 
\begin{equation} \label{lambd} \sum^l_{ i=0 } \frac{2c_i}{\sqrt{4\alpha^{2ip} e^{ 2 \lambda}+1}} = n \end{equation}
implies the following asymptotics of $ \lambda $ as $ n \to \infty  $,
\begin{equation}\label{lambdaas} \limsup_{ n\to \infty } \left| \lambda - \left(  - p \ln n + p \ln \left[ \frac{ \alpha - 1 }{ \ln \alpha } J ( p ) \right] \right) \right| = O( \alpha - 1)   , \end{equation} 
\[ J(p) = 2^{ 1 - 1/p } \int_0^\infty \frac{ du }{ \sqrt{ 1+u^{ 2p } } }. \]  

Step 11 - Plugging the asymptotics (\ref{lambdaas}) into (\ref{s4res}) we have,
\begin{equation}\label{estxi} n \left( s^{(4)}_n \right)^{ \frac 1{np} } = \frac{ \alpha - 1 }{ \ln \alpha } J( p ) \exp \left(  \frac\xi{np} + o(1) + O ( \alpha - 1 ) \right) , \end{equation}
\[ \xi \colon = \sum^l_{ i=0 } c_i \ln \frac{ \sqrt D_i + 1 }{ \sqrt D_i - 1 } . \]

Step 12 - For each $ \alpha > 1 $ the quantity $ \xi $ satisfies 
\[ n \frac{ 1 - \alpha^{ -2p }}{ 2( \alpha - 1 )} + o ( 1 ) \le \xi \le n \alpha  \frac{ \alpha^{ 2p } - 1}{ 2( \alpha - 1 )} + o ( 1 ) \]
as $ n \to \infty $. On plugging this in (\ref{estxi}) we find
\[ \limsup_{ n \to \infty } \left| n \left( s^{(4)}_n \right)^{ \frac 1{np} } - \frac{ \alpha - 1 }{ \ln \alpha } e J( p ) \right| =  O ( \alpha -1 ) . \]
Passing to the limit $ \alpha \to 1 $ we establish that 
\[ n \left( s (n) \right)^{ \frac 1{np} } \mathop{\longrightarrow}_{ n \to \infty } e  J( p ) .\]
It now remains to notice that the substitution $ v = \left( 1 + u^{ 2p }\right)^{ -1} $ in the definition of $ J ( p ) $ gives 
\[ J (p) = \frac 1p 2^{ - \frac 1p } B \left( \frac 1{2p} , \frac 12 - \frac 1{2p} \right) = \frac 1p B \left( \frac 1{2p} , 1- \frac 1{p} \right)  \]
on account of the identity $ \Gamma ( 2x ) = \frac{ 2^{ 2x-1}}{ \sqrt \pi } \Gamma ( x ) \Gamma \left( x + \frac 12 \right) $, see e. g.\cite{Artin}. The theorem is established. 
  
\subsection{Details}
 
  \bfseries
1)
\mdseries
Define $ \mR_ j \colon = \{ ( x_1 , \dots x_j ) \in {\mathbb N}^j \colon 1 \le x_1 \le x_2 < x_3 \le \ldots \prec x_{j} \} $, so
\begin{equation}  \delta_n \colon = s(n) - s_n^\prime =  \sum_{ \mR_n, \; x_n > T(n) }(x_1x_2\cdots x_{n})^{-p}  . \label{delta} \end{equation}
We have to show that 
  \begin{equation} \label{limsup}  \limsup \left( \frac 1n \ln \left( \frac{ \delta_n }{ s^\prime_n } \right) \right) \le  0 . \end{equation}
To this end, let us for $ l < n $ decompose the range of summation indices in (\ref{delta}) as follows, 
\[ \bigcup_{ j=1}^l \left\{ ( x_1 , \dots , x_n ) \in \mR_n \colon  \begin{array}{rcl} x_n &  > &  T ( n ) , \cr x_{ n-1 }&  > & T ( n-1  ) \cr & \cdots & \cr x_{ n - j + 1 } & > & T ( n -j + 1 ) \cr x_{ n-j } & \le & T ( n-j ) \end{array} \right\} \cup \left\{ ( x_1 , \dots , x_n ) \in \mR_n \colon \begin{array}{rcl}   x_n & > & T ( n ) , \cr x_{ n-1 } & > & T ( n-1  ) \cr & \cdots & \cr x_{ n - l } & > & T ( n -l ) \end{array}  \right\} , \]
and arrange the sum accordingly. This gives\footnote{The inequality in \eqref{estdelta}, as opposed to an equality,  is due to the fact that the $ x_j $'s in the last term on the right hand side ($ \xi_l s ( n - l ) $), written as a multiple sum, do not necessarily satisfy $ ( x_1 , \dots , x_n ) \in \mR_n $, say, $ x_{ n - l - 1} $ is not required to be less than $ x_{ n-l+1 } $.}  
\begin{eqnarray}\label{estdelta}  \delta_n \le  
\sum_{ j=1 }^{l-1} s^\prime_{ n-j } \xi_j + \xi_l s( n-l )   , \\ \nonumber  \xi_j \colon = \sum_{ \begin{array}{c}  x_n > T ( n ) , \cr x_{ n-1 } > T ( n-1  ) \cr \cdots \cr x_{ n - j + 1 } > T ( n -j + 1 ) \cr x_{ n-j+1 } \prec \dots \prec x_n \end{array} } (x_{ n-j+1} x_{ n-j+2 }\cdots x_{n})^{-p}  .
 \end{eqnarray} 
We estimate $ \xi_j $ by dropping all conditions on $ x_j $ except for the last two, $  x_{ n - j + 1 } > T ( n -j + 1 )$ and $ x_{ n-j+1 } \prec x_{ n-j+2 } \prec \dots \prec x_n $, and then repeatedly using the inequality \eqref{trivi}. Let' for definiteness, $ n $ be even, so $ x_n \ge x_{ n-1 }$, hence 
\begin{eqnarray*} \xi_j  \le \frac 1{p-1} \sum_{ \begin{array}{c}  x_{ n - j + 1 } > T ( n -j + 1 ) \cr x_{ n-j+1 } \prec \dots \le x_{ n-2 } <  x_{n-1} \end{array} } (x_{ n-j+1} x_{ n-j+2 }\cdots x_{n-2})^{-p}  ( x_{ n-1 } - 1 )^{ 1-2p } \le \\ \frac 1{2 (p-1)^2 } \sum_{ \begin{array}{c}  x_{ n - j + 1 } > T ( n -j + 1 ) \cr x_{ n-j+1 } \prec \dots < x_{ n-3 } \le  x_{n-2} \end{array} } (x_{ n-j+1} x_{ n-j+2 }\cdots x_{n-3})^{-p}  ( x_{ n-2 } - 1 )^{ 2-3p } \le \dots \\ 
\dots \le \frac 1{ j! ( p-1 )^j } \left( T ( n-j ) - \frac j2 \right)^{ j ( 1-p ) } . \end{eqnarray*}
Here the inequality holds whenever $ T ( n -j ) > n/2 $.  Let us plug this into (\ref{estdelta}) with $ l \sim n/2 $. For $ j \le l $ we have
\begin{equation} \xi_j \le \frac { C^{ j ( 1 -p)} }{ j! ( p-1 )^j } T \left( n-j \right)^{ j ( 1-p ) }  \label{xii} \end{equation}
with a constant $ C $ independent of $ A $. It follows that for $ j \le l $
 \begin{eqnarray*} s_{ n-j }^\prime \xi_j \le C^n s ( n-j ) \frac {  T \left( n-j \right)^{ j ( 1-p ) } }{ j! }\le  C^n \left( n-j \right)^{ Aj ( 1-p ) - p ( n-j ) } = C^n \left( n-j \right)^{ j [ A ( 1-p ) +p ] - np } . \end{eqnarray*}
 Here we have taken into account that $ s_{n-j}^\prime \le s(n-j) $, applied the upper estimate of (\ref{trivial}) to $ s ( n-j ) $, and droped $ j! $ in the denominator altogether. 
 
Notice now that the argument used in the the proof of the lower estimate in \eqref{trivial} shows in fact that $ s_n^\prime \ge C^n n^{ -np} $, for the value of $ ( x_1, x_2 , \dots , x_n ) $ used in that argument belongs to the domain of summation indices in the definition of $ s_n^\prime $. Applying this we get that   
 \[ \frac{ s_{ n-j }^\prime \xi_j }{ s_n^\prime } \le  C^n \left( n-j \right)^{ Aj ( 1-p ) - p ( n-j ) } n^{ np } \le C^n \left( n-j \right)^{ j [ A ( 1-p ) +p ] } . \]
In the last step we took into account that $ n-j \ge n/2 $ under our choice of $ l $.
 
 Let $ j \ge n \varepsilon_n $ with $ \varepsilon_n \downarrow 0 $ to be chosen later, and  let $ A > p/ ( p-1 ) $. In this case one can continue the inequality,
 \begin{equation} \label{largej}\frac{ s_{ n-j }^\prime \xi_j }{ s_n^\prime } \le C^n \left( \frac n2 \right)^{ n \varepsilon_n [ A ( 1-p ) + p ] } ,\end{equation}
 with a constant $ C $ still independent of $ A $.
 
 To deal with $ j < \varepsilon_n n $ we observe that 
 \begin{eqnarray*} s_n^\prime \ge \frac{ T( n )^{ 1-p } }{ p-1 } s^\prime_{ n-1 } \ge \dots \ge \frac 1{ ( p-1 )^j } \left( n ( n-1 ) \dots ( n-j+1 ) \right)^{ A ( 1-p) } s_{ n-j }^\prime   \end{eqnarray*}
 which upon substitution of (\ref{xii}) gives 
 \begin{eqnarray*}  \frac{ s_{ n-j }^\prime \xi_j }{ s_n^\prime } \le C \frac{ ( n-j )^{ A ( 1-p ) } }{j! \left( n ( n-1 ) \dots ( n-j+1 ) \right)^{ A ( 1-p)  } } = \\ 
\frac C{j!}  \exp \left( A ( p-1)  \sum_{ k=1}^j \log \left( 1 + \frac k{ n-j } \right) \right) \le C e^{ A ( p-1 ) \frac {j^2}{ n-j } } \le e^{ C n \varepsilon_n^2 } .
 \end{eqnarray*}
 Now, comparing the last displayed line with (\ref{largej}) and choosing $ \varepsilon_n = D / ( \log n ) $ with $ D $ large enough we find that  
 \[ \log_+ \frac{ s_{ n-j }^\prime \xi_j }{ s_n^\prime } =  o \left( \frac n{\ln n }  \right) \]
 uniformly in $ j $, and this implies the required estimate since (\ref{estdelta}) contains just $ O ( n ) $ terms.
 

  \bfseries
4)
\mdseries
Let us prove that the described map is a bijection between the set $ \{ ( x_1 , x_2 , \dots x_{ a_i }) \colon  x_1 \prec x_2 \ldots \prec x_{a_i} \} $ and the set of strictly increasing sequences $x'_1 < x'_2 <\ldots < x'_{a_i}\le w_i +c_i$. Firstly, the inverse transform is $ (x'_1, x'_2 \ldots ,x'_{a_i}) \mapsto ( x'_1-m_1, x'_2-m_2, \dots, x'_{a_i}-m_{a_i} )$ so the injectivity is obvious. Secondly, if a sequence $\{x_i\}$ satisfies the inequalities, then either $x_i<x_{i+1}$, or $x_i=x_{i+1}$ and the sign between $x_i$ and $x_{i+1}$ is $\le$, so $m_{i+1}>m_i$, and $x'_{i+1}>x'_i$. Thirdly, if a sequence $\{x'_i\}$ is strictly increasing, then either $m_i=m_{i+1}$,and $x_i<x_{i+1}$, or $m_i=m_{i+1}-1$, so $x_i\le x_{i+1}$ and and the sign between $x_i$ and $x_{i+1}$ is $\le$, so $\{x_i\}$ satisfy the inequalities.

  \bfseries
5)
\mdseries
We have
\[ (2T'(n))^2 {c_i+w_i \choose a_i} \ge {c_i+\frac{a_i}{2} \choose a_i } \ge \frac 1{(2T'(n))^2} {c_i+w_i \choose a_i }, \] 
and thus the ratio of values of $ H $ and $ H^\prime $ satisfies 
\[ (2T'(n))^{2l+2} \ge \frac{H( a_0, a_1 , \dots a_l) }{H^\prime ( a_0, a_1 , \dots a_l)} \ge \frac{1}{(2T'(n))^{2l+2}} . \] 
Since $ \left(2T'(n) \right)^{2l+2}=e^{O(\ln^2n)} $ this implies the assertion of this item.

 \bfseries
6)
\mdseries
The assertion follows from the fact that the number of terms in the sum does not exceed the number of subsets of $l+1$ integers lying between $0$ and $n$, which is $(n+1)^{l+1}=e^{O(\ln^2n)}=e^{o(n)}$.

 \bfseries
7)
\mdseries
First, $ s_n^{ (3)} \ge s^{ (2)}_n$ as maximum in the definition of $ s_n^{ (3) } $ is taken over a larger set. Second, let $(X_0, X_1, \ldots, X_l)$ be the point of maximum of $H'$, and let $Y_i=\lceil X_i\rceil$. We proceed in two steps,

(1) Replace all of $X_i$ by $Y_i$, and let $M=\sum Y_i - \sum X_i$. Notice that $M=O(\ln n )$. Let us write explicitly the expression for the maximal value of $ H^\prime $,   
\[ H'(X_0, X_1, \ldots , X_l) = \prod_{i=0}^l \frac{\alpha^{-iX_ip}\Gamma \left(\frac{X_i}{2}+c_i+1 \right)}{\Gamma \left( -\frac{X_i}{2}+c_i+1 \right)\Gamma(X_i+1)} . \] 
When replacing $X_i$ by $Y_i$, the value of $\alpha^{-iX_ip}$ will change by an $O(T'(n)^p )$ multiple, and the gamma-functions - by   $O(T'(n))$ multiples, so every factor will gain an at most an $e^{O(T'(n))}$ multiple, thus the whole product will change by no more than an $e^{O(\ln^2n)}=e^{o(n)}$ multiple; 

(2) Replace any nonzero $Y_i$ by $Y_i - 1$.  Again every operation will change the product in at most $e^{O(T'(n))}$ times, and the number of operations is  $O(\ln n )$.

After two these steps we will obtain a set of integers, $(X_0^\prime , X_1^\prime , \ldots , X_l^\prime )$ subject to $ \sum_{j=0}^l X_j^\prime = n $, for which  
\[
e^{o(n)} \ge \frac{H'(X'_0, X'_1, \ldots , X'_l)}{H'(X_0, X_1, \ldots , X_l)} \ge e^{-o(n)} , \] 
 and the assertion follows. 

\bfseries
8)
\mdseries
According to the Stirling formula, for every $x \le n$ we have 
\[ e^{O(-\ln(n))} \le \frac{e^{x\ln(x) - x}}{\Gamma(x+1)} \le e^{O(\ln(n))}, \]
with the implied constants independent of $x $ and $ n $. Thus, when
  $$\frac{\alpha^{-ia_ip}}{(c_i+\frac{a_i}{2}+1) B(a_i+1,c_i-\frac{a_i}{2}+1)}$$ is replaced
 by
 \[ \exp \left( - ia_ip\ln\alpha +\left(\frac{a_i}{2}+c_i\right)\ln\left(\frac{a_i}{2}+c_i\right)-a_i\ln a_i -\left(c_i-\frac{a_i}{2}\right)\ln(c_i-\frac{a_i}{2}) \right) \] 
 the product multiplies by an $e^{O(\ln(n))}$ factor, hence
\begin{eqnarray*} \prod^l_{i=0}\frac{\alpha^{-ia_ip}}{(c_i+\frac{a_i}{2}+1) B(a_i+1,c_i-\frac{a_i}{2}+1)} = \\ e^{o(n)}  \prod_{i=0}^l \exp\left( -\ln(\alpha)ia_ip+(\frac{a_i}{2}+c_i)\ln(\frac{a_i}{2}+c_i)-a_i\ln(a_i)-(c_i-\frac{a_i}{2})\ln(c_i-\frac{a_i}{2}) \right), \end{eqnarray*}
and the assertion follows. 

\bfseries
10)
\mdseries
First, notice that (see (\ref{lambd}))
\begin{equation} \label{lowerest} n = \sum^l \frac { 2 c_j }{ \sqrt D_j } \le ( \alpha -1 ) e^{ - \lambda } \sum^l  \alpha^{ i ( 1 -p ) } \le C ( \alpha ) e^{ - \lambda } . \end{equation} 
Then, observe that $ D_l \to \infty $ as $ n \to \infty $. Indeed, otherwise the sum $ \sum^l  \frac{ c_i}{\sqrt D_i } $ is estimated below by $ n^l \ge C n^A $ on a sequence of $ n $ for a non-zero $ C $.   Fix an $ M > 0 $ and let $ n $ be large enough  so as to there exists a $ j_M < l $ such that  $ D_j > M $ for $ j \ge j_M $, $ D_j \le M $ for $ j < j_M $. Clearly, $ l - j_M \to \infty $ as $ n \to \infty $ hence for all $ n $ large enough we are going to have,
\begin{eqnarray*}   \sum^l \frac{ 2 c_i }{ \sqrt D_i } \ge \sum^{ l }_{j_M} ( \dots ) \ge  e^{ - \lambda } ( \alpha - 1 ) \sum^l_{ j_M } \alpha^{ i ( 1-p)  } =   e^{ - \lambda } \frac{ \alpha - 1 }{ 1 - \alpha^{ 1-p}}   \alpha^{ j_M ( 1 -p ) }\left( 1  + o ( 1 )  \right) \ge \\   C_\alpha e^{ - \lambda - \lambda \frac{ 1-p }p } ( 1 + o ( 1 ))   . \end{eqnarray*}
Here $ C_\alpha $ is a constant in $ n $ uniformly separated from zero for small $ \alpha - 1 $.  Thus, there exists a positive $ C $ for all $ \alpha $ close enough to $ 1 $
\begin{equation}\label{lower} \limsup_{ n \to \infty } \frac{ e^{ -\lambda / p }}n < C .  \end{equation}
Proceeding, notice that 
\[ \sum^l \frac{ 2 \alpha^j}{ \sqrt{ 1 + 4 \alpha^{ 2jp} e^{ 2 \lambda }}} = \int_0^l \frac{ 2 \alpha^x dx}{  \sqrt{ 1 + 4 \alpha^{ 2px} }}+ O \left( e^{ -\lambda /p } \right)  . \]
To establish this it is enough to notice that the function $ h ( s ) = \frac s{ \sqrt{ s^{ 2p } + 1 }} $ has two regions of monotonicity over $ \mathbb R_+ $ separated by a point of maximum. It follows that the difference of the sum and the integral has absolute value not greater than trice the maximum of the integrand, and that the maximum is estmated above by $ C e^{ -\lambda / p } $ with the positive constant $ p $ depending on $ p $ only. On substituting, this and (\ref{lower}) give
\[ \limsup_{ n\to \infty } \left| 1 - ( \alpha - 1 ) \frac 1n  \int_0^l \frac{ 2 \alpha^x dx}{  \sqrt{ 1 + 4 \alpha^{ 2px} }} \right| = O \left(  \alpha - 1  \right) . \]
On the other hand, 

\[  \int_0^l \frac{ 2 \alpha^x dx}{  \sqrt{ 1 + 4 \alpha^{ 2px} }} = e^{ - \lambda/p }  \frac 1{ \ln \alpha } \int_{ e^{\lambda/p} 2^{ 1/p}}^{\alpha^l 2^{ 1/p } e^{ \lambda / p } } \frac{ 2^{ 1 - 1/p } du}{  \sqrt{ 1 + u^{ 2p} }} . \]
The lower limit of integration vanishes in the limit $ n \to \infty $ by  (\ref{lowerest}), the upper one goes to infinity by (\ref{lower}) andthe fact that $  \alpha^l \asymp n^A $ with an $ A > 1 $, hence the integral in  the right hand side goes to $ J ( p ) $ as $ n \to \infty $. Thus,  
\[ \limsup_{ n\to \infty } \left| 1 - \frac{\alpha - 1 } { \ln \alpha }\frac{ e^{ - \lambda/p } }n  ( J(p) + o ( 1 ) ) \right| = O \left(  \alpha - 1  \right) , \]
which is (\ref{lambdaas}). 

\bfseries
12)
\mdseries
Applying the summation by parts we find
\[ \xi = \left( \sum_{ m=0 }^l c_m \right)\ln \frac{ \sqrt{ D_ l} + 1 }{  \sqrt{ D_l} - 1 } + \sum^l_{ j=0} \left( \sum_{ m=0 }^j c_m \right) \left[ \ln \frac{ \sqrt{ D_{ i-1 }} + 1 }{  \sqrt{ D_{ i-1 }} - 1 } - \ln \frac{ \sqrt{ D_ i} + 1 }{  \sqrt{ D_i} - 1 } \right] . \]
By the definition of $ c_i $, $ \sum_{m=0}^j  c_m = \lfloor \alpha^{ j+1 } \rfloor $ for any $ j $, hence the first term in the right hand side is $ O ( \alpha^{l ( 1-p )} e^{ -\lambda } ) = O ( n^{ A ( 1-p ) +p } ) $ and hence vanishes as $ n \to \infty $ with our choice of $ A $.  The second term is estimated by applying the mean value theorem to the difference in the square brackets, which gives
 \[ [ \dots ] = \alpha^{ 2 ( i-i_*)p} \frac { 1 - \alpha^{ -2p }}{ \sqrt { 1 + 4 \alpha^{ 2 i_* p }e^{ 2\lambda } }} ,\;  i_* \in ( i-1 , i ) . \] 
 Replacing the denominator with $ \sqrt{ D_i } $, then $\sqrt{ D_{ i-1  }} $, and taking into account that $ \sum^l \frac { 2 c_i }{ \sqrt {D_i}} = n $, we obtain
 \[ n \frac{ 1 - \alpha^{ -2p }}{ 2( \alpha - 1 )} \le \sum^l \alpha^j \left[ \ln \frac{ \sqrt{ D_{ i-1 }} + 1 }{  \sqrt{ D_{ i-1 }} - 1 } - \ln \frac{ \sqrt{ D_ i} + 1 }{  \sqrt{ D_i} - 1 } \right] \le n \alpha  \frac{ \alpha^{ 2p } - 1}{ 2( \alpha - 1 )} , \]
 and the assertion follows.

\bfseries
\begin{center}
Overview
\end{center}
\mdseries
The above proof uses the special growth function $x^{p}$ in the model problem.  Still, some points in the above argument admit generalisation.  We enumerate them according to the steps in the proof.
\bfseries
2)
\mdseries
Consider  $$s_n=  \sum_{1 \le x_1 \le x_2 < x_3\le \ldots \prec x_{n}\le T(n)}(F(x_1)F(x_2)\cdots F(x_{n}))$$
for some functions $F$ and $T$, when $T(n)$ is positive integer and $F$ is positive, $k_n=(s_n)^{\frac{1}{np}}$, $p>0$ is some real number. Our aim is to calculate $k=\lim_{n \to \infty}k_n$. Consider for every $\alpha>1$ the sequence of integers $1=c_0(\alpha)\le c_1(\alpha) \le \ldots$ and the function $P_n(\alpha)=F(c_j(\alpha))$, where $j$ is the maximal integer for which $c_j(\alpha)\le n$. Let  $$s'_n(\alpha)=  \sum_{1 \le x_1 \le x_2 < x_3\le \ldots \prec x_{n}\le T(n)}(P_{x_1}(\alpha)P_{x_2}(\alpha)\cdots P_{x_n}(\alpha))$$ and $k'_n(\alpha)=(s'_n(\alpha))^{\frac{1}{np}}$
Then if for every $\epsilon>0$ exist such $\alpha_0$ that for every $1<\alpha<\alpha_0$ if $P_n(\alpha)=P_m(\alpha)$ then $\frac{F(n)}{F(m)}<1+\epsilon$ we can apply the same argument like in section 2) and conclude that if  exist sequences $A_n(\alpha)$, $B_n(\alpha)$, $C_n(\alpha)$ and function $G(\alpha)$ with following properties:
 \bfseries
1)
\mdseries
$A_n(\alpha)<k'_n(\alpha)<B_n(\alpha)$
 \bfseries
2)
\mdseries
$\lim_{\alpha \to 1}\lim_{n \to \infty}B_n(\alpha)-A_n(\alpha)=0$
 \bfseries
3)
\mdseries
$A_n(\alpha)<C_n(\alpha)<B_n(\alpha)$
 \bfseries
4)
\mdseries
$\lim_{n \to \infty}C_n(\alpha)=G(\alpha)$
 \bfseries
5)
\mdseries
$\lim_{\alpha \to 1}G(\alpha)=k'$ for some real $0<k'<\infty$
Then $k=k'$.

\bfseries
3)
\mdseries
This section does not depend on the specific growth function at all.
 Let us fix some $\alpha>1$ and $n$. Let us rename $s'_n(\alpha)$ by $s_n$, $P_n(\alpha)$ by $P(n)$ and $c_n(\alpha)$ by $c_n$. In the sum \[ s_n=  \sum_{1 \le x_1 \le x_2 < x_3 \le \ldots \prec x_{n}<T(n)}(P(x_1)P(x_2)\cdots P(x_{n}))^{-pn}\] each $P(x_i)$ is equal to one of the numbers $1=P(c_0),P(c_1), \ldots, P(T(n))$, and $l$ is maximum integer for which $c_l(\alpha)\le T(n)$. Let $d_i=c_{i+1}-c_i$ be the number of such $y$ that $P(y)=P(c_i)$. Then $$s_n=\sum_{all sequences a_0 , a_1, \ldots ,a_l with \sum{a_i}=n}H(a_0, a_1, \ldots , a_n)P(c_0)^{a_0}P(c_1)^{a_1}\cdots P(c_l)^{a_l}$$ when $H(a_0, a_1, \ldots, a_l)$ - number of ways to choose numbers $0<x_1\le x_2\le <x_3 \le x_4 <\ldots \prec x_n<T(n)$, such that for every $i$ numbers of $j$: $P(x_j)=\alpha^{i}$ is equal to $a_i$, because for this 
$0<x_1\le x_2\le x_3< x_4 \le \ldots \prec x_n<T(n)$ $P(x_1)P(x_2)\cdots P(x_n)=P(c_0)^{a_0}P(c_1)^{a_1}\cdots  P(c_l)^{a_l}$. 

\bfseries
4)
\mdseries
The assertions of this section again do not depend on the specific function.
We can apply same argument as above and conclude, that $H(a_0, a_1, \ldots, a_l=\prod^l_{i=0}C^{a_i}_{d'_i+w_i}$, where $d'_i=d_i=c_{i+1}-c_i$ if $i<l$, $d'_l=T(n)-c_l$, and  $w_i=\lfloor \frac{a_i-1}{2}\rfloor$ or $w_i=\lceil \frac{a_i-1}{2}\rceil$, depending on parities of $\{a_0, a_1, \ldots, a_l\}$. $w_i=\lceil \frac{a_i-1}{2}\rceil$ if $a_i$ is even and $\sum_{j=0}^ia_j$ is even, else $w_i=\lfloor \frac{a_i-1}{2}\rfloor$. 

\bfseries
\begin{center}
Acknowledgements
\end{center}
\mdseries
I should express my gratitude to Roman Romanov, who introduced this problem for me and  helped me with solving the problem and writing the column. The work was supported by the RSF Grant RSF 17-11-01064.

\end{document}